\documentclass[twocolumn]{autart}    

\usepackage{graphicx}          

\usepackage{amsmath}
\usepackage{amsfonts}
\usepackage{amssymb}
\usepackage{mathrsfs}
\usepackage{enumitem}
\usepackage{scrextend}

\newtheorem{theorem}{Theorem}[section]

\newtheorem{lemma}[theorem]{Lemma}

\theoremstyle{definition}

\newtheorem{remark}{Remark}

\newcommand{\enurom}[1]
{\begin{enumerate}[label=\rm{(\roman*)}]
		#1
\end{enumerate}}

\newcommand{\dd}[0]
{\text{d}}

\newcommand{\enualp}[1]
{\begin{enumerate}[label=\rm{(\alph*)}]
		#1
\end{enumerate}}

\newcommand{\ccd}[0]
{(.)}

\newcommand{\ttnn}[1]
{\textnormal{#1}}

\newcommand{\modulo}[1]
{\left|#1\right|}

\newcommand{\rif}[1]
{(\ref{#1})}

\newcommand{\oo}[0]
{\infty}

\newcommand{\norm}[1]
{\left\| #1 \right\|}

\newcommand{\ra}
{\rightarrow}

\newcommand{\PS}[2]
{\langle\,#1,#2\rangle}

\newcommand{\ps}[2]
{#1 \cdot #2}

\newcommand{\rs}[0]
{(s)}

\newcommand{\mineq}[1]
{\leq #1}

\newcommand{\mageq}[1]
{\geq #1}

\newcommand{\equazioneref}[2]
{
		\begin{eqnarray}\label{#1}
		\begin{split}
		#2
		\end{split}
		\end{eqnarray}
}

\newcommand{\sistemaref}[2]
{
		\begin{eqnarray}
		\begin{cases}\label{#1}
		#2
		\end{cases}
		\end{eqnarray}
}

\newcommand{\sistemanoref}[1]
{
		\begin{eqnarray*}
			\begin{cases}
				#1
			\end{cases}
		\end{eqnarray*}
}

\newcommand{\eee}[1]
{\begin{eqnarray*}\begin{split}
				#1
\end{split}\end{eqnarray*}}

\newcommand{\graffe}[1]
{\left\{ #1 \right\}}

\newcommand{\tonde}[1]
{\left ( #1 \right )}

\newcommand{\ccal}[1]
{\mathcal{#1}}

\newcommand{\scr}[1]
{\mathscr{#1}}

\newcommand{\bb}[0]
{\textit I\!}

\newcommand{\ldue}[3]
{\scr L^2(#1,#2;\bb R^{#3})}

\newcommand{\acc}[0]
{\` }

\begin{document}

\begin{frontmatter}

\title{Lipschitz regularity of controls and inversion mapping for a class of smooth extremization problems} 


\author[Basco]{Vincenzo Basco}\ead{vincenzo.basco@thalesaleniaspace.com}    

\address[Basco]{Thales Alenia Space, Via Saccomuro, 24 Roma Italy}  


\begin{abstract}                          
In the contest of optimal control problems, regularity results for optima are known when addressing fiber-strictly convex Lagrangian. For infinite time horizons, or for settings with infinite dimensional dynamics, the equivalence between minima/maxima and extremals could break down. Commonly, this is due to a loss of convexity/concavity  of the cost functional or to a presence of state constraints, in which further controllability assumptions are needed. For many science applications, this a trend is not required, as in energy saving problems. In the present paper, we deal with the set of a functional's extremals subject to end-point restrictions. We consider an affine control system and a cost functional associated to an autonomous Lagrangian. The dynamics is smooth, satisfying the Lie bracket condition, and the functional is assumed merely Fréchet differentiable. Here we provide a regularity result for controls in the context of constrained extremization problems, under weaker conditions on Lagrangian than the not met classical ones. More precisely, we show a  characterization for the Lipschitz regularity of controls associated with the extremal trajectories steering two fixed points, assuming the absence of singular controls. As main application, we construct a locally Lipschitz inversion mapping from the ambient space to the set of constrained extremals.  
\end{abstract}

\end{frontmatter}

\section{Introduction}
Classical control theory takes into account minimization -- or maximization -- problems of functionals along trajectories given by a dynamical system (cfr. \cite{bardi2008optimal}, \cite{montgomery2006tour}, \cite{vinter2010optimal}). These optimization problems are generally well defined in the presence of convex -- or concave -- data or for non-negative Lagrangians. For many science applications, this a trend is not required for energy saving problems, where non-convex energy functionals and the presence of singularities are commonly an indication of instabilities, so that the direct methods in classical optimization cannot be applied (see e.g. \cite{conti2017analytical},\cite{maso2008vanishing},\cite{lin2006simulations}). In this case, it is required that a system evolves along a trajectory that makes the functional action stationary, i.e. its first variation is zero. Since any critical point of a smooth convex -- or concave -- function is a minimum -- or a maximum,  the extremization problems reduce so to classical optimization problems when the functional meets these requirements and the constraints are linear. Hence, the seek of  {extremals} -- or stationary points --  rather than minima -- or maxima -- reveal to be necessary, the concept of extremization is a natural generalization of optimization to the non-convex case. This formulation of control problems in terms of stationary optimization has recently proved in \cite{bascodowerexploiting2021} to be very useful for the investigation of boundary constrained problems.

The point of view of extremization problems appears useful in engineering applications, where the effects of small perturbations are not negligible. As in the case of classical optimization problems, a fundamental role for regularity investigations of feedback controls is played by the absence of singularity (see Section 2 below). Singular points may occur in  control problems: the importance of singular paths -- namely associated to a singular control -- of affine control system was investigated since a long time in the in calculus of variations (cfr. \cite{bellaiche1996sub},\cite{bonnard2003singular},\cite{calin_chang_2009}). Singular paths are candidates as minimizers for time optimal control problems (cfr. \cite{bardi2008optimal},\cite{rifford2014sub}) and they are strictly related to the singularities of the {end-point map}. It is a well known fact that for control unconstrained optimization problems associated to running cost functional, the solutions of the {abnormal maximum principle} projects onto the set of singular trajectories. However, such set and the solutions of the normal maximum principle  does not trivially intersect generally: this could happen when considering degenerate affine control systems with end-point constraints (e.g., the geodesics of the manifold spanned by the three-dimensional Martinet distribution).

The optimal synthesis of control problems is quite challenging while continuity of the cost functional, sufficient and necessary conditions, and sensitivity relations in absence of singular paths are known (cfr. \cite{agrachev2002strong},\cite{bardi2008optimal},\cite{vinter2010optimal},\cite{bascocannfrank2019semisubRiem},\cite{cannarsa2008semiconcavity}) and duality connections, involving the Legendre transform, are investigated (see \cite{ekeland1977legendre}). The classical theory on the regularity of controls in the absence of singularities has been exhaustively invested for Lagrangian two times differentiable and strictly convex in the control variable (\cite{bellaiche1996sub} and the literature therein). Nevertheless, the absence of singularity does not fit for the optimal synthesis in many investigations in applied science. To deal with it, in the case of optimization problems, well known and succesfull results have been developed by using continuation arguments for indirect algorithms methods or convergently probability-one homotopy techniques (cfr. \cite{watson2001theory},\cite{allgower2012numerical}).

Although the assumption of absence of singular trajectories could be strong, results to ensure regularity of feedback controls under less restrictive assumptions on the Lagrangian are needed when addressing extremization under end-point constraints. We emphasize, in optimization context, singular time optimal paths imply the failure of the Lipschitz continuity of feedback controls. This  can be see by \cite[Proposition 5.3- (ii)]{bascocannfrank2019semisubRiem} which implies, along a singular time optimal paths, the horizontal proximal gradient of time minimal function containing a non-zero vector. A necessary condition for the Lipschitz continuity of feedback controls at a given end-point is the triviality of horizontal proximal gradient of time minimal function at such point. So, without making additional assumptions -- such as the absence of singular minimizing controls like in \cite{bascocannfrank2019semisubRiem} -- one cannot hope for a Lipschitz regularity.

In this paper we relax the classical assumptions, requiring at most the well definition of the problem.
The main concern is to provide a regularity characterization for controls associated with extremal trajectories of a given cost functional. More exactly, we show that the set of controls related to end-point constrained extremals is equi-Lipschitz continuous, in absence of singular paths. The trajectories are assumed to solve an affine dynamics and the Lagrangian is possibly unbounded and not strictly convex in the fiber.
We provide further an application of this result constructing a Lipschitz inversion mapping for constrained extremals.

The outline of this paper is as follows. In Section 2 we give some basic definitions, notations, and the general background.  Section 3 indicates the problem statement and the main assumptions.  In Section 4 we state the main result of this paper showing a regularity characterization for controls associated with constrained extremal trajectories. The Section 5 is devoted to an application case.

\section{Background}\label{section_2}

For any $u,w\in \bb R^n$ we denote by $| u|$ and $\ps{ u}{ w}$ the Euclidean norm of $u$ and the scalar product between $u$ and $w$, respectively. $\bb B_\delta (x)$ is the closed ball in $\bb R^m$ centered at $x$ of radius $\delta>0$, and $\bb B:=\bb B_1(0)$, $\bb B_\delta:=\bb B_\delta (0)$. For any $p\mageq 1$, $\scr L^p(a,b;\bb R^m)$ denotes the space of all Lebesgue measurable functions on $[a,b]$ that are $p$-integrable, endowed with the norm $\norm{u}_{p}^p=\int_a^b |u\rs|^pds$. The Hilbert space $\ldue{a}{b}{m}$ is equipped with the scalar product defined by $\PS{u}{w}=\int_a^b \ps{u\rs}{w\rs}ds$. For any subset $A\subset \ldue{a}{b}{m}$, we denote by $\ttnn{cl}_w A$ the weak closure of $A$.
Consider $(V,\norm{.}_V)$ and $(B,\norm{.}_B)$ be two Banach spaces and let $\scr V\subset V$ be an open subset. Consider a mapping \(\Phi: \scr V \rightarrow B\). If $\Phi$ is linear, we denote by $\Phi^\star$ the adjoint (or transpose) operator. 
We denote by $\dd \Phi$ the {Fr\'{e}chet derivative} of $\Phi$ (whenever it exists). High order Fr\'{e}chet derivatives (whenever they exist) of $\Phi$ at $u$ shall written $\dd ^j \Phi(u)$ or $\dd ^j_u \Phi(u)$. If $\dd^j\Phi$ is a continuous map on $\scr V$ (wrt the operator norm topology) we say that $\Phi$ is of class $C^j$.
Let $\scr F$ be a family of smooth vector fields on $\bb R^n$. We denote by $Lie_{\scr F}(x)$ the \textit{Lie algebra} of vector fields generated by $\scr F$ at $x\in \bb R^n$ constructed as follows. Let $Lie^1_{\scr F}(x)$ the set defined by $ \ttnn{span }\{X(x)\,|\, X\in \scr F$
and define recursively, for all $i\in \bb N^{+}$, the set $Lie^{i+1}_{\scr F}(x)$ defined by the spanning of $Lie^i_{\scr F}(x)\cup  \{[X,Y](x)\,|\, X\in Lie^1_{\scr F}(x),\, Y\in Lie^i_{\scr F}(x)\}$
where $[X,Y](x):=\dd Y(x)(X(x))-\dd X(x)(Y(x))$ is the Lie bracket. Then the Lie algebra is defined by $Lie_{\scr F}(x)= \bigcup_{i\in \bb N} Lie^{i+1}_{\scr F}(x)$.

\subsection{Extremization problems}

Let $X$ be an Hilbert space and consider a family of subset $\{\scr V(x)\}_x\subset X \times X \times$ $\bb{R}$. We denote by  
 \eee{
\delta \scr V(x) =0
 }
the problem of determining all couples $(u, z) \in X \times \bb{R}$ such that
\equazioneref{ext_problem}{
(u, 0, z) \in \scr V(x)\text {. }
}
The set of all $z \in X$ such that there exists $u \in X$, called extremal, satisfying \rif{ext_problem} is denoted by 
$
\{ \delta \scr V(x) =0\}
$
and we refer to it as the extremization (or staticization) problem associated to $\scr V(x)$. It occurs e.g. when $\scr V(x)$ is the Lagrangian sub-manifold associated with some function $\Phi: X \rightarrow \bb R$ that is $C^1$ (in the Fréchet sense)
$
\scr V(x)=\left\{\left(u, \dd \Phi (u), \Phi(u)\right) \mid u \in X\right\}
$
and $\dd \Phi$ stands for the Fréchet derivative. In that case \rif{ext_problem} becomes
$
\dd \Phi (u)=0$ and $  z=\Phi(u),
$
so that $\{ \delta \scr V(x) =0\}$ is simply the unconstrained smooth extremization problem of determining the critical points and values of $\Phi$. 
We consider the particular case when $\{\scr V(x)\}_{x\in \bb R^n}$ is the family of manifolds
\eee{
&\{ \,(u, \lambda_0\dd \Phi (u)-\sum_{j=1}^{n} \lambda_{j} \dd_u \ccal C_{j} (u;x), \Phi(u) ) \text{ such that} \\
&\;\ccal C(u;x)=0, (\lambda_0,\lambda)\in \bb{R}^{1+n}-\{0\}\,\}
}
where $\ccal C(x)=(\ccal C_j(.;x))_{j}$ and $f,\ccal C_{j}(.;x):X\ra \bb R$ are $C^{1}$ functions for any $x$. If  $\{\dd_u \ccal C_{j} (u;x), 1 \mineq j \mineq n\}$ is a set of linearly independent vectors for any
\eee{
u \in S(x):=\{\ccal C(.;x)=0\}
}
it follows that $u\in \{ \delta \scr V(x) =0\}$ iff $u$ is a critical point of $\left.\Phi\right|_{S(x)}$, the restriction of $\Phi$ to $S(x)$. In that case, we write $\delta \scr V (x)=0$ as
\equazioneref{ext_phi_S}{
\delta_{S(x)} \Phi=0
}
and call it a constrained smooth extremization problem. Whenever the set $\{\delta_{S(x)} \Phi=0\}$ is non-empty, we define the set-valued function
\eee{
x\mapsto \langle\delta\rangle_{S(x)} \Phi : =\{\Phi(u) \,|\, u\in \{\delta_{S(x)} \Phi=0\}\}.
}
If $S(x)$ is the whole space, we shall write $\langle\delta\rangle_u \Phi $ in place of $\langle\delta\rangle_{S(x)} \Phi $.

\subsection{Lack of convexity: a motivational example}

It is of great interest to consider extremization problems rather that minima/maxima for a large family of functional as those expressed by energy-driven (or saving) systems. Such class of functionals are represented by
\eee{
&\Phi(u)=  \int_{I} V(\omega, u)+W(\omega,\nabla u) d \omega, \quad I\subset \bb R^n
}
where $V$ and $W$ are continuous function and possibly non-convex. As example consider the family of Ginzburg–Landau type energy functionals
\eee{
\Phi_{GL}(u)=\int_{I}|u^2-\alpha|^{2}+q^{2}| |\nabla u|^2-\beta|d\omega,\quad \alpha,\beta,q \mageq 0
}
arising in many frameworks, including convective pattern formation and magnetic thin films. It is straingforward to show that end-point constraints reduce the functional $\Phi_{GL}$ to be generally non-convex. To see that, consider the one-dimentional case $I=[-1,1]$ with parameters $\alpha=\beta=q=1$, control law
\eee{
u'=v, \quad v \in \scr L^2(-1,1;\bb R)
}
and subject to end-points constraints
\eee{
&u(-1)=u(1)=0.
}
We have, $\Phi_{GL}(0)=4$ and $\Phi_{GL}(\xi_{\pm})=16/15$ where $\xi_{\pm}(\omega)=\pm (1-|\omega|)$.

\subsection{Singular points}
Consider $X_1,...,X_m$ smooth vector fields on $\bb R^n$ and the solution $\xi(.)$  associated with $u \in \scr L^2(0,T;\bb R^m)$ of the affine control system
\equazioneref{sistema_affine}{
&\xi'\rs=\sum_{i=1}^{m}u_i\rs X_i(\xi\rs),\; {s\in[0,T]}\ttnn{ a.e.}\\
&\xi(0)=x_0\in \bb R^n
}
with the end-point constraints
\equazioneref{C_equazione}{
\ccal C(u;x)= \xi (T)-x.
}
A point $u\in S(x)$ is said to be a singular point for $S(x)$ if
\eee{
\dd_u \ccal C(u;x) \text{ is not surjective.}
}
We denote by $S'(x)$ the set of all singular points of $S(x)$. Singular points may occurs in the set of extremization problem $\delta_{S(x)} \Phi=0$, e.g., in sub-Riemannian structures.
The importance of \textit{singular paths}  (namely associated to a singular control, i.e. a singular point of $S(x)$) of the above affine control system
was investigated since a long time in the in calculus of variations (cfr. \cite{bellaiche1996sub},\cite{bonnard2003singular},\cite{calin_chang_2009}). The singular paths are candidates as minimizers for time optimal control problems (cfr. \cite{bardi2008optimal},\cite{rifford2014sub}  and the literature therein)
and they are strictly related to the singularities of the  {end-point map} (cfr. Section 2). Indeed, for control problems minimizing a running cost functional associated to $L$, bounded from below in the state, such trajectories are solutions of the {abnormal ($\lambda_0=0$, see below) maximum principle}  corresponding to the associated Hamiltonian
\eee{
\sup_{u\in \bb R^m} (\sum_{i=1}^{m}u_i\ps{p}{X_i(x)}+\lambda_0 L(x,u)).
}
When $m=n$, as known, extremals are not abnormal, i.e. $\lambda_0\neq 0$. The optimal synthesis of such control problems is quite challenging although continuity of the cost functional, necessary conditions, and sensitivity relations in absence of singular paths are known (cfr. \cite{agrachev2002strong},\cite{bardi2008optimal},\cite{vinter2010optimal}) and duality connections, involving the Legendre transform, are investigated (see \cite{ekeland1977legendre}).

Nevertheless, in the case $m<n$, there exist extremals that are normal but singular (cfr. \cite{bellaiche1996sub}). Indeed,
in this setting, there are privileged and prohibited paths, and the singular trajectories turn out to be the singularities of the set of curves satisfying the constraint \rif{C_equazione}.
%
In particular, when final end-point conditions are imposed, the sub-Riemannian manifold generated by the vector fields $X_i$ may have points for which there are not curves steering them. To deal with it, H\"{o}rmander in \cite{hormander1967hypoelliptic} introduced a condition on the Lie algebra associated to the distribution spanned by $X_i$'s (see Section 2).


\section{Problem statement and main assumptions}
Let $T>0$ and $X_1,...,X_m$ a set of vector fields on $\bb R^n$. We consider  the following
\begin{assum}\label{ass_1}  $X_i$'s are smooth ($C^\infty$ or analytic) bounded vector fields, linearly independent any $x\in \bb R^n$, and
$Lie_{\graffe{X_1,...,X_m}}(x)=\bb R^n$ for all $x\in \bb R^n$.
\end{assum}

\begin{remark}
\noindent The mapping
$
	x\mapsto \ttnn{span} \graffe{X_1(x),...,X_m(x)}
$
is called a \textit{smooth distribution} of rank $m$ on $\bb R^n$ and the condition on the Lie algebra in Assumption \ref{ass_1} is called {H\"{o}rmander's condition}. From the Chow-Rashevsky Theorem (see  \cite{chow1940systeme},\cite{rashevsky1938connecting}), for any $x,y\in \bb R^n$ there exists an absolutely continuous arc $\xi:[0,T]\ra \bb R^n$, with square integrable derivative, such that 
\eee{
&\xi'\rs \in \ttnn{span} \graffe{X_i(\xi\rs)}_i\quad s\in [0,T] \ttnn{ a.e.}\\
&\xi(0)=x, \,  \xi(T)=y.
}
Furthermore, for any arc $\xi\ccd$ steering $x\in \bb R^n$ to $y\in \bb R^n$ in time $T$, there exists a unique $u_\xi\in \scr L^2(0,T;\bb R^m)$ satisfying
\eee{
\xi'(s)=\sum_{i=1}^m (u_{\xi}\rs)_i X_i(\xi\rs), \;s\in [0,T] \ttnn{ a.e.}
}
\end{remark}
For any $x_0\in \bb R^n$ consider $\scr V^{T,x_0}\subset \scr L^2(0,T; \bb R^m)$ be such that for every control $u\in \scr V^{T,x_0}$ the associated solution $\xi_{x_0,u}\ccd$ of the affine system \rif{sistema_affine}, with starting point $\xi_{x_0,u}(0)=x_0$, is well defined on $[0,T]$. Since $X_i$'s are smooth, then $\scr V^{T,x_0}$ can be chosen to be open (see  \cite{bellaiche1996sub,montgomery2006tour}). The \textit{end-point map} $E_{T,x_0}:\scr V^{T,x_0}\ra \bb R^n$ is defined by
\eee{
E_{T,x_0}(u)= \xi_{x_0,u}(T).
}
Notice that, under assumptions Assumption 1, for any $x_0\in \bb R^n$ we can choose
\eee{
\scr V^{T,x_0}=\ldue{0}{T}{m}.
}

 Let $L:\bb R^n\times \bb R^m\ra \bb R$ be a two time continuously differentiable Lagrangian. In the following, we fix $T>0$, $x_0,x\in \bb R^n$, and denote by $\Phi:\ldue{0}{T}{m}\ra \bb R$ the functional  defined by
\eee{
	\Phi(u) =\int_0^T L(\xi_{x_0,u}\rs,u\rs)\,ds
}
and by $\ccal C:\ldue{0}{T}{m}\ra \bb R^n$ the map
\eee{
\ccal C(u)=E_{T,x_0}(u)-x.
}
Throught this paper, in order to carrying out our analysis, we consider the basic conditions
\equazioneref{condizione_non_triviale}{
\Phi\in C^1 \quad \&\quad\{\delta_{S(x)} \Phi=0\} \neq \emptyset
}
and the further assumptions
\begin{assum}\label{ass_2}
for any $x\in \bb R^n$, $\dd_u L(x,.)$ is a $C^1$-diffeomorphism, i.e. bijective with continuously differentiable inverse.
\end{assum}

\begin{assum}\label{ass_3}
$\exists \theta ,\,\psi,\,\varphi:\bb R^+\ra\bb R^+$ continuous satisfying 
\eee{
&\liminf_{|u|\ra \oo} \frac{\theta(|u|)}{|u|^2}>0\\
&\limsup_{|x|\ra \oo} \frac{\psi(|x|)}{|x|^2}<+\oo
}
such that
\eee{
	&L(x,u)\mageq \theta(|u|)-\psi(|x|),\, \forall x\in \bb R^n\, \forall u\in \bb R^m \\
	&|\dd_xL(x,u)| \leqslant  \varphi(r)\left(|u|^{2}+1\right),\, \forall x \in \bb B_r\, \forall u \in \bb{R}^{m}.
}
\end{assum}
\begin{assum}\label{ass_4}
$S(x)' \cap  \, \{\delta_{S(x)} \Phi=0\} =\emptyset$.
\end{assum}
It follows some remarks.
\begin{remark}
\
\enualp{\label{remark2}
\item It is well known that, under the Assumption 1, $E_{T,x_0}$ turn out to be of class $C^1$ and $\dd E_{T,x_0}$ is weakly continuous (see \cite{bellaiche1996sub}). The {H\"{o}rmander's condition} on the vector fields $X_i$'s needs to be imposed to consider controllability for end-point constraints. If $\ttnn{span } \{X_i(x)\}_i$ is the whole space, e.g. as in mass-spring configurations systems, then it is satisfied.
\item From the regularity of $(\Phi,\ccal C)$, Remark \ref{remark2}-(a), and the definition of extremal given in the previous section, it straightforward to show that the set $\{\delta_{S(x)} \Phi=0\}$ is weakly closed wrt the $\scr L^2$ topology. Furthermore, \rif{condizione_non_triviale}, Assumption 4, and Lagrange's Multiplier Theorem leads to consider such set be non-empty.
\item Applying the Inverse Mapping Theorem, Assumption 2 is satisfied for the class of strictly convex Lagrangian in the control, i.e. $\dd_u^2 L>0$. In such a setting, the condition \rif{condizione_non_triviale} is naturally met for minimizing control problems. Moreover, under the Assumption 2 and 4, and by using the Pontryagin's maximum principle, the regularity characterization of minimizing controls are known (cfr. \cite{bonnard2003singular},\cite{calin_chang_2009}). 
\item From Assumption 4, we would like to remark that the set of Lagrange multipliers $(\lambda_0,\lambda)\neq (0,0)$ reduces to consider $\lambda_0\neq 0$.
}
\end{remark}

\section{Regularity of constrained extremals}
	   Next, we state and prove the main result of this paper under the Assumptions 1-4. We provide a characterization of regularity of controls associated with constrained extremals in the space of uniformly Lipschitz continuous functions, under the assumption of absence of singular controls.

\begin{thm}[Regularity characterization]\label{regularity_control_result} Consider Assumptions \ref{ass_1}-\ref{ass_4}.
	Then the following statements are equivalent:
	\enurom{
\item $\Phi|_{\{\delta_{S(x)} \Phi=0\}}<+\infty$.
	
		\item there exists  \(K>0\) such that
		\eee{
			|u\rs|\mineq K\quad \forall s\in [0,T]\,\forall u  \in\{\delta_{S(x)} \Phi=0\}.
		}
		
		\item  there exists  \( K>0\) such that 
		\eee{
			&|u(T)|\mineq  K,\quad |u(s)-u( t)|\mineq  K|s-t|\\
			&\forall s, t\in [0,T]\, \forall u  \in\{\delta_{S(x)} \Phi=0\}.
		}
		
	}

\end{thm}

\subsection{Some preparatory results}
In the following we define the set-valued Hamiltonian \(H: \bb{R}^{n} \times \bb{R}^{n}  \multimap \bb{R}\) by
\eee{
	H(x, p)=\langle\delta\rangle_u (\,\sum_{i=1}^m u_i\ps{p}{ X_i(x)}- L(x, u)\,  )
}
and, for all $x,p\in \bb R^n$ and $u\in \bb R^m$, we put
\eee{
	&h(x, p, u):= \sum_{i=1}^mu_i\ps{p}{ X_i(x)}- L(x, u)\\
	&Z(x, p):=(\ps{p}{X_1(x)},...,\ps{p}{X_m(x)})^{\star}.
}

\noindent Furthermore, for any $u\in \{\delta_{S(x)} \Phi\}$, we define the (non-empty) set
\eee{
	\Lambda (u;x) :=\{   \lambda \in \bb R^n\,|\, \dd \Phi(u)=\lambda \dd_u \ccal C(u;x) \}.
}

\begin{lem}\label{lemma_1}
	The set-valued Hamiltonian is single-valued and there exists a locally Lipschitz continuous function ${w}:\bb{R}^{n} \times \bb{R}^{m} \rightarrow \bb{R}^{m}$, with ${w}(x,.)$ a $C^1$-diffeomorphism for any $x\in \bb R^n$, satisfying:
	\enurom{
	
	\item ${w}(x,.)=(\dd_u L(x,.))^{-1}$ for all $x\in \bb R^n.$
	
	\item \(H(x,p)={w}(x, Z(x, p))\cdot Z(x, p)-L(x, {w}(x, Z(x, p)))\) for all $x,p\in \bb R^n$.
}
\end{lem}
\begin{pf}
	From Assumption 2 and the regularity of $L$, there exists a locally Lipschitz map ${w}:\bb{R}^{n} \times \bb{R}^{m} \rightarrow \bb{R}^{m}$ such that ${w}(x,.)$ is a $C^1$-diffeomorphism for all $x\in \bb R^n$ and satisfying the statement (i). Moreover, for all \((x, p, u)\in \bb{R}^{n} \times \bb{R}^{n} \times \bb{R}^{m} \)
	\sistemaref{differenziale_h_x_p_u}{
		\dd_{x} h(x, p, u)=\sum_{i=1}^{m} u_{i} \dd X_{i}(x)^{\star} p-\dd_{x} L(x, u) \\ 
		{\dd_{p} h(x, p, u)=\sum_{i=1}^{m} u_{i} X_{i}(x)} \\ 
		{\dd_{u_{i}} h(x, p, u)=\ps{p}{X_i(x)}-\dd_{u_{i}} L(x, u)},\; \forall i=1,...,m
	}
	which implies that, for any \((x, p) \in \bb{R}^{n} \times \bb{R}^{n}\), $\dd_{u} h(x, p, u)=0 $ if and only if
	\eee{
		\dd_{u} L(x, u)=(\ps{p}{X_1(x)},...,\ps{p}{X_m(x)})^{\star}.
	}
	So, for any \((x, p, u) \in \bb{R}^{n} \times \bb{R}^{n} \times \bb{R}^{m}\), we have that $\dd_{u} h(x, p, u)=0$ if and only if $	u={w}(x, Z(x, p))$. Hence, the Hmiltonian is single-valued and for all $x,p\in \bb R^n$
	\equazioneref{u_controllo_massimo_in_H}{
		&H(x,p)\\
		&=h(x, p, {w}(x, Z(x, p))\\
		&= {w}(x, Z(x, p))\cdot Z(x, p)-L(x, {w}(x, Z(x, p)))
	}
and the statement (ii) follows.
\end{pf}

\begin{prop}\label{lemma_p_co_stato_u_uguale_u_gotico}
	Let $u\in \{\delta_{S(x)} \Phi\}$.
		
	Then, for any $\lambda\in \Lambda (u;x) $, there exists an absolutely continuous arc \(p=p_u:[0, T] \rightarrow \bb{R}^{n}\) solving the Hamiltonian system
	\equazioneref{Ham_system}{		
		\begin{cases}
				\xi'\rs=\dd_{p} H(\xi\rs, p\rs)\\
			-p'\rs=\dd_{x} H(\xi\rs, p\rs)
		\end{cases}
	\quad s\in[0,T]\ttnn{ a.e.}
}
	with final condition
	\equazioneref{cond_finale_p}{
		p(T)=\lambda
	}
which satisfies the relation
\equazioneref{B_star_p_uguale_Z}{
	Z(\xi\rs,p\rs)=\dd_uL(\xi\rs,p\rs)\quad \forall s \in[0, T].
}
In particular, the control $u\ccd$ satisfies
	\equazioneref{controllo_u_ottimo}{
		u(s)=  {w}  \left(\xi(s), Z\left(\xi(s), p(s)\right)\right)\quad \forall s \in[0, T].
	}
	
\end{prop}
\begin{pf}
	Let \(v  \in \scr L^{2}\left(0, T ; \bb{R}^{m}\right)\) and consider $\lambda\in \Lambda (u;x) $. We have
	\eee{
	&\PS{\dd \Phi( u)}{v}\\
	&=\int_{0}^{T} \dd_xL(\xi\rs,u\rs)\cdot \dd E_{s,x_0}(u)(v)   \\
	&\qquad\qquad\qquad\qquad +  \dd_u L\left(\xi(s), u(s)\right)\cdot v(s) \, d s.
}
	Let us denote by $\Psi\ccd$ the fundamental solution of
	\sistemanoref{
		\Psi'\rs=A\rs \Psi\rs & s\in [0,T]\ttnn{ a.e.} \\
		\Psi(0)=I
	}	
	where for all $s\mageq 0$
	\eee{
		&A\rs:=\sum_{i=1}^m u_i\rs  \dd X_i(\xi\rs)\\
		&B\rs:=(X_1(\xi\rs),...,X_m(\xi\rs)).
	}
	Then, by the known expression of $\dd E_{s,x_0}(u)$ (cfr. \cite[Proposition 1.8]{rifford2014sub}),
	\eee{
		&\int_{0}^{T}   \dd_xL(\xi\rs,u\rs)\cdot \dd E_{s,x_0}(u)(v)   d s\\
		&=\int_{0}^{T} (\dd_xL(\xi\rs,u\rs)\cdot  \int_{0}^{s} \Psi(s) \Psi(r)^{-1} B(r) v(r) d r ) d s\\
		&=\int_{0}^{T} \int_{0}^{s}   \dd_xL(\xi\rs,u\rs)\cdot \Psi(s) \Psi(r)^{-1} B(r) v(r)  d r d s.
	}
	Using the Fubini theorem, we have
	\eee{
		&\int_{0}^{T} \int_{0}^{s}   \dd_xL(\xi\rs,u\rs)\cdot \Psi(s) \Psi(r)^{-1} B(r) v(r)  d r d s\\
		&=\int_{0}^{T} \int_{r}^{T} \dd_xL(\xi\rs,u\rs)\cdot \Psi(s) \Psi(r)^{-1} B(r) v(r)  d s d r.
	}
	Hence,
	\eee{
		&\int_{0}^{T}   \dd_xL(\xi\rs,u\rs)\cdot \dd E_{s,x_0}(u)(v)   d s\\
		&=\int_{0}^{T}  ( \int_{s}^{T}\left(\Psi(r) \Psi(s)^{-1} B(s)\right)^{\star} \\
		&\qquad \qquad \qquad \cdot \dd_xL(\xi(r),u(r)) d r\cdot v(s) )   d s\\
		&=\int_{0}^{T}  (   B(s)^{\star}\left(\Psi(s)^{-1}\right)^{\star}\\
		 &\qquad \qquad \qquad \cdot \int_{s}^{T} \Psi(r)^{\star} \dd_xL(\xi(r),u(r)) d r  \cdot   v(s) )  d s.
	}	
	Furthermore, for all \(v  \in \scr L^{2}\left(0, T ; \bb{R}^{m}\right)\)
	\eee{
	&\lambda\cdot \dd E_{T,x_0}(u)(v) \\
	&=  \lambda\cdot \int_{0}^{T} \Psi(T) \Psi(s)^{-1} B(s) v(s) d s   \\ 
	&=\int_{0}^{T}   \lambda\cdot \Psi(T) \Psi(s)^{-1} B(s) v(s) d s  \\ 
	&=\int_{0}^{T} B(s)^{\star}\left(\Psi(s)^{-1}\right)^{\star} \Psi(T)^{\star} \lambda\cdot v(s)   d s. 
	}
	Let us set for every $s\in [0,T]$
	\eee{
	&p(s):=\left(\Psi(s)^{-1}\right)^{\star} \Psi(T)^{\star} \lambda -\left(\Psi(s)^{-1}\right)^{\star}\\
	&\qquad \qquad \qquad \cdot \int_{s}^{T} \Psi(r)^{\star} \dd_xL(\xi(r),u(r)) d r.
	}
	Then, since $v$ is arbitrary, it follows that
	\[
	B(s)^{\star} p(s)=\dd_{u} L\left(\xi(s), u(s)\right) \quad \forall s \in[0, T]
	\]
	and so \rif{B_star_p_uguale_Z} holds. Now, the relation in \rif{controllo_u_ottimo} follows from \rif{B_star_p_uguale_Z} and Lemma \ref{lemma_1}-(i). Also,
	\[
	- {p}'(s)=A(s)^{\star} p(s)-\dd_x L(\xi\rs,u\rs) \quad s \in[0, T]\text { a.e.} 
	\]
	So, recalling \rif{u_controllo_massimo_in_H} and from \rif{controllo_u_ottimo}, we have for a.e. $s\in [0,T]$
\sistemanoref{
	\dd_{x} H(\xi\rs, p\rs)=\dd_{x} h\,[s] \\ 
	\dd_{p} H(\xi\rs, p\rs)=\dd_{p} h\,[s],
}
where we put $[s]:=(\xi\rs, p\rs, {w}(\xi\rs, Z(\xi\rs, p\rs)))$ for all $s\in [0,T]$. Then, recalling \rif{differenziale_h_x_p_u}, the pair $(\xi\ccd,p\ccd)$ solves the Hamiltonian system \rif{Ham_system} with final condition \rif{cond_finale_p}.
\end{pf}

\subsection{Proof of Theorem \ref{regularity_control_result}}

The implication (ii)$\Longrightarrow$(i) immeditely follows from the equibounded property of trajectories and from the continuity of
the Lagrangian. We need only to show (ii)$\Longrightarrow$(iii) and (i)$\Longrightarrow$(ii).
	
%
	
	(ii)$\Longrightarrow$(iii): From the embedding $ \scr L^2(0,T;\bb R^m)\subset \scr L^1(0,T;\bb R^m)$ and applying the Gronwall Lemma,  we can choose \(\hat K>0\) such that
	\equazioneref{bound_on_u_in_l_uno_e_due_e_su_xi}{
		&\norm{u}_{1}+\norm{u}_2\mineq \hat  K,\quad \xi_{x_0,u}(.)\subset \bb B_{\hat K}\\
		 &\;\forall u \in \{\delta_{S(x)} \Phi=0\}.
	}	
	First, we claim that the family of co-states $p_u \ttnn{ solving } \rif{Ham_system}$
	\equazioneref{co_states_p_u_set}{
		&\{\,p_u:[0,T]\ra \bb R^n\,|\,  p_u(T)=\lambda\in \Lambda (u;x) ,\\
		&\;u\in \{\delta_{S(x)} \Phi=0 \}\,\}
	}
	is equibounded. We have that
	\eee{
		\exists c>0\,:\,\left|\lambda \right| \leqslant c \quad \forall \lambda \in \Lambda (u;x) \,\forall  u \in \{\delta_{S(x)} \Phi=0\}.
	}
	Otherwise, suppose, by contradiction, that there exist two sequences \(\graffe{u_k}_k\subset  \{\delta_{S(x)} \Phi=0\}\) and $\lambda_k\in \Lambda (u_k;x) $ such that \(\lim_k\left|\lambda_k\right|= +\infty\). We may assume that $\lambda_k\neq 0$ for all $k\in \bb N$. Thus, for every \(k\in \bb N\) and $v  \in \scr L^{2}(0, T; \bb{R}^{m})$
	\equazioneref{p_k_d_E_uguale_d_J}{
		\frac{\lambda_k}{\left|\lambda_k\right|}\cdot \dd_u \ccal C\left(u_{k};x\right)(v )   = \frac{1}{\left|\lambda_k\right|}   \PS{\dd \Phi\left(u_{k} \right)}{v}.
	}	
	Passing to subsequences and keeping the same notation, we may assume that  there exist \( \lambda \in \bb{R}^{n}\) and $u\in \ttnn{cl}_w\,\{\delta_{S(x)} \Phi=0\}$, with  \(|\lambda|=1 \), such that $ 	\frac{\lambda_{k}}{\left|\lambda_{k}\right|}\ra \lambda$ and $u_k\ra u$ weakly. Recalling Remark \ref{remark2} and passing to the limit as \(k \rightarrow \infty \) in \rif{p_k_d_E_uguale_d_J} we have
	\eee{
		\lambda\cdot \dd_u \ccal C (u ;x) = 0
	}
	and, from our assumptions, a contradiction follows. Now, from the smoothness of the vector fields $X_i$, let	\(M>0\) be such that \(|\dd X_{i}(x)^{\star}\xi| \leqslant M|\xi|\) for all $\xi \in \bb R^n$, \(x \in \bb B_{\hat K}\), and all \(i=1, \ldots, m\). Applying Proposition \ref{lemma_p_co_stato_u_uguale_u_gotico}, we have that for any $u\in \{\delta_{S(x)} \Phi=0\}$ and any $\lambda \in \Lambda (u;x) $ there is an absolutely continuous arc \(p_{u} :[0, T] \rightarrow \bb{R}^{n}\) such that  for a.e. $s \in[0, T] $
	\sistemanoref{
		{p}_{u}'(s)  = \sum_{i=1}^{m}u_i(s)  \dd X_{i}\left(\xi_{x_0,u}(s)\right)^{\star} p_{u}(s)\\
		\qquad \qquad\qquad\qquad\qquad -\dd_{x} L\left(\xi_{x_0,u}(s), u(s)\right)\\
		p_{u}(T)=\lambda.
	}
	 By Assumption 3, we can find a constant $\tilde K=\varphi({\hat K})>0$ such that  for every $u\in \{\delta_{S(x)} \Phi=0\}$ and for a.e. $  s \in[0, T] $
	\eee{
	\left| {p}_{u}'(s)\right| \leqslant M \sqrt{m}|u(s)| \left|p_{u}(s)\right|+\tilde K\left(|u(s)|^2 +1\right).
	}
	Thus, for all $u\in \{\delta_{S(x)} \Phi=0\}$ and all $  s \in[0, T] $
	\eee{
		&\left|p_{u}(s)\right| \\
		&\mineq \left|p_{u}(T)\right|+\int_{s}^{T} \tilde K\left(|u(r)|^2+1\right) d r\\
		&\qquad +\int_{s}^{T} M \sqrt{m}|u(r)| \left|p_{u}(r)\right| d r \\ 
		&\mineq (c+\tilde K(\hat K^2+T))+\int_{s}^{T} M \sqrt{m}|u(r)|\left|p_{u}(r)\right| d r .
	}
	Applying Gronwall's Lemma and recalling \rif{bound_on_u_in_l_uno_e_due_e_su_xi}, it follows the claim \rif{co_states_p_u_set}. 
%
	Thus, from the claim and the Gronwall's Lemma again, we have that there exists $ R\mageq K$ such that for any $u\in \{\delta_{S(x)} \Phi=0\}$ and any co-states $p_u\ccd$ in \rif{co_states_p_u_set}
	\equazioneref{bound}{
		\xi_{x_0,u}\ccd\subset \bb B_{ R},\quad p_u\ccd \subset \bb B_{ R}.
	}
Hence, using \rif{controllo_u_ottimo} and the locally Lipschitz continuity of $\mathfrak u(.,.)$, we get (iii).

	(i)$\Longrightarrow$(ii): Fix $u\in \{\delta_{S(x)} \Phi=0\}$ and, from Assumpation 4, let $C>0$ and $\alpha>0,\hat  \alpha>0$ such that $\theta(r)\mageq \alpha r^2$ and $\psi(r)\mineq \hat \alpha r^2$ for all $r> C$. For all $t\in [0,T]$, we have
		\equazioneref{stima_di_partenza}{
		\int_0^t\modulo{u(s)}^2ds&=\int_{[0,t]\cap \graffe{s:\modulo{u(s)}> C}}\modulo{u(s)}^2\,ds\\
		&\qquad+\int_{[0,t]\cap \graffe{s:\modulo{u(s)}\mineq C}}\modulo{u(s)}^2\,ds\\
		&\mineq\int_{[0,t]\cap \graffe{s:\modulo{u(s)}> C}}\modulo{u(s)}^2\,ds+C^2T.
	}
Now, from (i) and Assumption 3, there exists $\hat M>0$, not depending on $u\ccd$, such that
\equazioneref{stima_per_gronwall_su_norma_u}{
&\int_{[0,t]\cap \graffe{s:\modulo{u(s)}> C}}\modulo{u(s)}^2\,ds\\
&\mineq \alpha^{-1} \int_{[0,t]\cap \graffe{s:\modulo{u(s)}> C}}\theta (\modulo{u(s)})\,ds\\
& \mineq \alpha^{-1} \int_0^t \theta (\modulo{u(s)})\,ds\\
&\mineq \alpha^{-1} \int_0^t L(\xi_{x_0,u}(s),u(s))\,ds\\
&\qquad+\alpha^{-1} \int_0^t \psi(|\xi_{x_0,u}\rs|)\,ds\\
&\mineq \alpha^{-1}\hat M+\alpha^{-1} \int_0^t \psi(|\xi_{x_0,u}\rs|)\,ds.
}
Putting $M=\max_{i} \sup_{x\in \bb R^n}|X_i(x)|<\infty$, we have for all $s\in [0,T]$
\eee{
|\xi_{x_0,u}\rs|&\mineq |x_0|+\int_0^s \sum_{i=1}^m |X_i(\xi_{x_0,u}(\tau))||u_i(\tau)|d\tau\\
& \mineq |x_0|+mM \int_0^s \sum_{i=1}^m |u_i(\tau)|d\tau.
}
Since $\psi(|\xi_{x_0,u}\rs|)\mineq \hat \alpha |\xi_{x_0,u}\rs|^2$ for all $s\in [0,T]$ and using that $\ldue{0}{T}{m}\subset \scr L^1(0,T;\bb R^m)$, we have that there exists a constant $\hat c>1$, not depending on $u\ccd$, such that for all $s\in [0,T]$ 
\eee{
&\psi(|\xi_{x_0,u}\rs|)\\
&\mineq \hat c\;(\int_0^s|u(\tau)|^2d\tau+|x_0| ({\int_0^s|u(\tau)|^2d\tau})^{1/2}+|x_0|^2).
}
Hence, for all $t\in [0,T]$
\equazioneref{2_stima_per_gronwall_su_norma_u}{
&\int_0^t \psi(|\xi_{x_0,u}\rs|)ds\\
&\mineq 2\hat c T|x_0|^2+2\hat c \int_0^t(\int_0^s|u(\tau)|^2d\tau)ds.
}
Now,
from \rif{stima_di_partenza}-\rif{2_stima_per_gronwall_su_norma_u}, 
the Gronwall's Lemma, and since all the involved constants does not depends on $u\in \{\delta_{S(x)} \Phi=0\}$, we conclude that there exists $ k>0$ such that $\norm{u}^2\mineq k$ for all $u\in \{\delta_{S(x)} \Phi=0\}$. Consequently, applying again Gronwall's Lemma, we have that all trajectories \(\xi_{x_0,u}(.)\)
	associated with controls \(u\in \{\delta_{S(x)} \Phi=0\}\) are uniformly bounded, that is, there exists $\tilde K>0$ such that for	every \(u \in \{\delta_{S(x)} \Phi=0\}\)
	\equazioneref{bound_on_traj}{
	\xi_{x_0,u}(.) \subset \bb B_{\tilde K}.
}	
	From Proposition \ref{lemma_p_co_stato_u_uguale_u_gotico} and replacing $\tilde K>0$ with a suitable greater constant, we have that for any $u\in \{\delta_{S(x)} \Phi=0\}$ and any co-states $p_u\ccd$ in \rif{co_states_p_u_set} satisfies
	\equazioneref{bound_on_p}{
	p_u\ccd \subset \bb B_{\tilde K}.
}
	Hence, from the continuity of $H(.,.)$, we get
	\eee{
	 N:=\sup_{(s,u)\in [0,T]\times \{\delta_{S(x)} \Phi=0\}}\;|H(\xi_{x_0,u}\rs,p_u\rs)|<\oo.
}
Now, notice that, from Proposition \ref{lemma_p_co_stato_u_uguale_u_gotico}, for all $u\in \{\delta_{S(x)} \Phi=0\}$ and all $s\in [0,T]$ we have
\eee{
&{w}  \left(\xi_{x_0,u}\rs, Z\left(\xi_{x_0,u}\rs, p(s)\right)\right)\cdot Z(\xi_{x_0,u}\rs,p\rs)\\
&=u\rs\cdot \dd_uL(\xi_{x_0,u}\rs,p\rs).
}
Hence, applying Lemma \ref{lemma_1}-(ii),  for all $u\in \{\delta_{S(x)} \Phi=0\}$ and all $s\in [0,T]$
\equazioneref{eq_finale}{
&H(\xi_{x_0,u}\rs,p_u\rs)\\
&=\ps{u(s)}{\dd_uL(\xi_{x_0,u}\rs,u(s))}-L(\xi_{x_0,u}\rs,u(s))).
}
Now,
%
%
	fix $R>0$ and consider \(\hat{N}>0\) be such that
	$
		|L(x, u)| \leqslant \hat{N}$ for all $x\in \bb B_R, u\in \bb B.
	$
	Since $L$ is smooth, for any $x\in \bb B_R$ and any $r>0$ the function $L(x,.)+c(x,r)|\,.\,|^2$ is convex on $\bb B_r$ for a suitable $c(x,r)\mageq 0$ (cfr. \cite{bardi2008optimal}). Moreover, fix $r=r(R)>0$ such that 
	\equazioneref{condizione_finale_star}{
	\theta(|u|)>\psi(|x|)+2\hat N+N,\; \forall |u|>r\, \forall x\in \bb B_R.
	}	
	Consider \(x \in \bb B_R, u \in \bb {R}^{m} \backslash\{0\},\) and put \(v:=u /|u|\). 
	Denote by $g=g_v:\bb R\ra \bb R$ the function $g(s):=\ps{\dd_u L(x,sv)}{sv}-L(x,sv)$. Then, from Assumption 2, the function $g$ is continuously differentiable. So, $g(s)=g(1)+\int_1^s g'(\tau)d\tau$ for all $s\mageq 1$. Now, for any $\tau\mageq 1$, we have $g'(\tau)=\tau\ps{\dd^2_uL(x,\tau v)v}{v}$. Hence, for all $s\mageq 1$, applying the Fundamental Theorem of Calculus we have
	\equazioneref{stella1}{
	&\ps{\dd_u L(x,sv)}{sv}-L(x,sv)\\
		&\mageq \ps{\dd_u L(x,v)}{v}-L(x,v)+\int_1^s \ps{\dd^2_u L(x,\tau v)v}{v}\; d\tau\\
		&=d_uL(x,v)\cdot v -L(x,v)+\tau d_uL(x,\tau v)\cdot v|^{s=\tau}_{1=\tau}\\
		&\qquad -\int_1^s \dd_u^2 L(x,\tau v) v\cdot v d\tau\\
		&=\dd_uL(x,sv)\cdot (sv-v)+d_uL(x,v)-L(x,v).
	}
	 Hence, since $|v|=1$, for any $\sigma>r$
	\equazioneref{stella2}{
	&\dd_uL(x,sv)\cdot (sv-v)\\
	&\mageq L(x,sv)-L(x,v)+c(x,\sigma)(s-1),\, \forall s\in [1,\sigma].
	}
	From \rif{stella1} and \rif{stella2}, we have
	\eee{
	&\ps{\dd_u L(x,sv)}{sv}-L(x,sv)\\
	&\mageq L(x,sv)+c(x,\sigma)(s-1)-2\hat N,\, \forall s\in [1,\sigma].
	}
	So, for any $\sigma>r$ and putting $s=|u|$
	\eee{
	&\ps{\dd_u L(x,u)}{u}-L(x,u)\\
	&\mageq \theta(|u|)-\psi(|x|)+c(x,\sigma)(|u|-1)-2\hat N,\, \forall |u|\in [1,\sigma].
	}
	Notice that $\sup_{x\in \bb B_R}c(x,\sigma)<+\infty$ for any $\sigma>r$. Summing up and recalling \rif{condizione_finale_star}, we can conclude that there exists
a constant \(K >0\) (depending only on $\tilde K$) satisfying
	\eee{
	\ps{\dd_u L(x,u)}{u}-L(x,u)>{N},\quad \forall x \in \bb B_{\tilde K}\, \forall |u|>K.
}
	By \rif{eq_finale}, we deduce that $|u(s)| \leqslant K$ for all $s \in[0, T]$ and all $u \in \{\delta_{S(x)} \Phi=0\}$. The proof of Theorem \ref{regularity_control_result} is now complete.


\section{Application to Existence of Lipschitz Inversion Mappings}
As main application of the Characterization Theorem \ref{regularity_control_result}, in this section we provide an existence results of Lipschitz inversion mapping. In the following we denote by $\scr L_{k}^T$ the class of all $\bb R^m$-valued $k$-Lipschitz continuous functions on the time interval $[0,T]$.

The proof of the following two lemmata are omitted since straightforward and well known.
\begin{lemma}\label{lemmino}
	Let $X$ be a separable normed space and $\Phi:X\ra \bb R^n$ be a linear, continuous, and surjective operator. Consider $\graffe{x_i}_i$ dense in $X$.	Then there exist linearly independent vectors $x_1,...,x_n$ such that $\Phi:W\ra \bb R^n$ is an isomorphism, where $W=\ttnn{span}\,\graffe{x_1,...,x_n}$.
\end{lemma}
\begin{lemma}[\cite{bascocannfrank2019semisubRiem},\cite{rifford2014sub}]\label{mappa_diff_3}
	Let $\tau>0$. Then the map
	\[[0,\tau] \times  L^2(0,\tau;\bb R^m) \ni (s,u)\mapsto \ttnn dE_{s,x_0}(u)\]
	is continuous.
\end{lemma}

\begin{thm}[Inversion mapping]\label{inversione}
%
	If $\Phi|_{\{\delta_{S(x)} \Phi=0\}}<+\infty$, then there exist $k> 0$, $r>0$, and $\ell\mageq 0$ such that for any 
	$
	(t,u)\in [0,T]\times \{\delta_{S(x)} \Phi=0\}
	$
	we can find a map
	\eee{
	\mathfrak{u}_{t,u}: B_{r}(t) \times B_{r}(x_{x_0,u}(t))\ra \{k\text{-Lip. func. on $[0,T]$}\}
	}
	satisfying:
	\enurom{
		\item $\mathfrak{u}_{t,u}\in C^1$.
		
		\item $E_{s,x_0}( \mathfrak{u}_{t,u}(s,\beta))=\beta$
		for all $ (s,\beta)\in \bb B_{r}(t) \times \bb B_{r}(x_{x_0,u}(t))$.
		
		\item $\ttnn d\mathfrak{u}_{t,u}$ is ${\ell}$-Lipschitz.
	}
\end{thm}

\begin{pf}
	Let $(t_0, u_0) \in [0,T]\times \{\delta_{S(x)} \Phi=0\}$. We know that
	$
	\ttnn dE_{t_0,x_0}( u_0)( \cdot)
	$
	is surjective on $\scr L^2(0,T;\bb R^m)$. Let $\scr V\subset C^1(0,T;\bb R^m)$ be a countable subset such that $\overline {{\ttnn{span} \scr V}}=\scr L^2(0,T;\bb R^m)$. Applying Lemma \ref{lemmino}, there exist $n$ linearly independent vectors $\graffe{v_1^{0},...,v_n^{0}}\subset \scr V$ such that the map
	\eee{
	(\alpha_i)_i\mapsto  \sum_{i=1}^n \alpha_i \ttnn dE_{t_0,x_0}( u_0)( v_i^0)
	}
	is an isomorphism on $\bb R^n$.
	Define for any $(t,u)\in (0,T) \times \scr L^2(0,T;\bb R^m)$ the map $\varphi^0_{t ,u}:\bb R^n\ra \bb R^n$ by	$\varphi^0_{t,u}((\alpha_i)_i)= \sum_{i=1}^n \alpha_i \ttnn dE_{t,x_0}( u)( v_i^0)$. By Lemma \ref{mappa_diff_3}, there exist $\varrho_0>0,\,\mu_0>0$ such that for any $(t, u)\in \ccal J_0:=[t_0-\rho_0,t_0+\rho_0] \times B_{\scr L^2}(u_0,\varrho_0)$ the map $\scr E_{t,u}^0:(0,T) \times \bb R^n\ra (0,T) \times \bb R^n$, defined by
	\eee{
	\scr E_{t ,u}^0(s ,(\alpha_i)_i)= (s ,E_{s,x_0}( u+\sum_{i=1}^m \alpha_iv_i^0))
	}
	 satisfies for all $(t, u)\in \ccal J_0$
	\eee{
	|\det\tonde{\ttnn d \scr E_{t, u}^0(t, 0)}|=\modulo{\det \varphi^0_{t ,u}}\mageq \mu_0
	}
	where $B_{\scr L^2}(u,r)$ stands for the closed ball in $\scr L^2$ centered at $u$ of radius $r$. Now, from the Ascoli-Arzel\acc a theorem, the set $[0,T]\times \{\delta_{S(x)} \Phi=0\}$ is compact. Then there exists $N\in \bb N^+$ such that, for all $j=1,...,N$, we can find $\rho_j>0,\,\mu_j>0$, $(t_j,u_j)\in [0,T]\times \{\delta_{S(x)} \Phi=0\}$, and linearly independent $\{v_1^j,...,v_n^j\}\subset \scr V$, such that
	\eee{
	&[0,T]\times \{\delta_{S(x)} \Phi=0\}\\
	& \subset \bigcup_{j=1,...,N}\, [t_j-\varrho_j,t_j+\varrho_j]  \times B_{\scr L^2}(u_j,\varrho_j)\\
	&=: \bigcup_{j=1,...,N} \ccal J_j.
	}
	Defining for any $(t, u)\in \ccal J_j$ the maps $\scr E_{t, u}^j:(0,T) \times \bb R^n\ra (0,T) \times \bb R^n$ and $\varphi^j_{t, u}:\bb R^n\ra \bb R^n$ by
	\eee{
	&\scr E_{t, u}^j(s, (\alpha_i)_i)= (s, E_{s,x_0}(u+\sum_{i=1}^m \alpha_iv_i^j))\\
	&\varphi^j_{t, u}((\alpha_i)_i)= \sum_{i=1}^n \alpha_i \ttnn dE_{t,x_0}( u)( v_i^j)
	}
	we deduce that for all $(t, u)\in \ccal J_j$ and $j=1,...,N$
	\equazioneref{stima}{
		&|\det\tonde{\ttnn d \scr E_{t, u}^j(t, 0)}|\\
		&=\modulo{\det \varphi^j_{t, u}}
		\mageq \mu_j\mageq \min\graffe{\mu_1,...,\mu_N}>0.
	}
%
%
%
	Applying the Inverse Mapping Theorem to the map $\scr E_{t ,u}^j$ and using a compactness argument, we conclude that for each $j$ there exists $r_j>0$ such that the set $  V_j(t, u):=(t-r_j,t+r_j) \times \bb  B_{r_j}(E_t( u))$ is isomorph to $ ({\scr E_{t, u}^j})^{-1}\tonde{  V_j(t, u)}$ for any $(t, u)\in \ccal J_j$. Put $r= \min\graffe{r_1,...,r_N}$ and define for any $(t, u)\in \ccal J_j$
	\eee{
		&\mathfrak{u}_{t, u}(t' ,\beta)= u+\sum_{i=1}^m \alpha_i(t' ,\beta)v_i^j\\
		&\forall (t' ,\beta)\in   V_j(t ,u)
	}
	where $(\scr E_{t, u}^j)^{-1}(t' ,\beta)=:(t' ,\alpha(t' ,\beta))$. Notice that, since the coefficients ${\alpha_i}$ are bounded by a suitable constant $M\mageq 0$ and ${v_i^j}\in C^1(0,T;\bb R^m)$, there exists a constant $k\mageq 0$ such that $\mathfrak{u}_{t ,u}$ take values in $\scr L_{k}^T$. Hence, $(i)$ and $(ii)$ follow. Moreover, from \rif{stima} and the regularity of the end-point map, there exists a constant $\ell\mageq 0$, depending only on $k$, such that $\ttnn d\mathfrak{u}_{t ,u}$ is $\ell$-Lipschitz for all $(t ,u)\in [0,T]\times \{\delta_{S(x)} \Phi=0\}$. So, we get $(iii)$.
\end{pf}

\section{Conclusions}
In this paper, we have addressed a class of extremization problems wrt end-point constraints qualifications. We assumed the functional generating the Lagrange sub-manifold merely Fréchet differentiable and the associated Lagrangian possibly unbounded and not strictly convex in the fiber.  In this paper, with the use of recently investigated techniques for extremization problems, we show that the controls associated with constrained extremals enjoy of Lipschitz regularity whenever singular trajectory are absent. {H\"{o}rmander's} condition on the affine dynamics cannot be skipped in order to the well-posedness of the problem. As main implication, we construct a locally Lipschitz inversion mapping from the ambient space to the set of constrained extremals, that turns out to be useful for the well definition of value functions associated to extremization problems.

\bibliographystyle{plain}
\bibliography{BIBLIO_VB_crit_traj}




\end{document}